\documentstyle[12pt]{amsart}
\newtheorem{theorem}{Theorem}[section]

\newtheorem{lemma}[theorem]{Lemma}

\numberwithin{equation}{section}

\begin{document}

\title[Energy identity for the maps from a surface with tension field bounded in $L^p$]
{Energy identity for the maps from a surface with tension field
bounded in $L^p$}

\author{Li Jiayu;\hspace{1cm} Zhu Xiangrong}

\address{Li Jiayu, Department of Mathematics, University of Science and Technology of China, Hefei 230026, P. R. China}
\email{lijia@@amss.ac.cn}

\address{Zhu Xiangrong, Department of Mathematics, Zhejiang Normal University, Jinhua 321004, P. R. China}
\email{zxr@@zjnu.cn}

\thanks {The research was supported by NSFC (11071236, 11101372).}

\begin{abstract}
Let $M$ be a closed Riemannian surface and $u_n$ a sequence of
maps from $M$ to Riemannian manifold $N$ satisfying
$$\sup_n(\|\nabla u_n\|_{L^2(M)}+\|\tau(u_n)\|_{L^p(M)})\leq \Lambda$$
for some $p>1$, where $\tau(u_n)$ is the tension field of the
mapping $u_n$.

For the general target manifold $N$, if $p\geq \frac 65$, we prove
the energy identity and neckless during blowing up.

\end{abstract}
\maketitle

\section{Introduction}
\hspace{6mm} Let $(M,g)$ be a closed Riemannian manifold and $(N,h)$ be a Riemannian manifold without boundary. For a mapping $u$
from $M$ to $N$ in $W^{1,2}(M,N)$, the energy density of $u$ is defined by
$$e(u)=\frac 12|du|^2=\textrm{Trace}_g u^*h$$
where $u^*h$ is the pull-back of the metric tensor $h$.

The energy of the mapping $u$ is defined as
$$E(u)=\int_M e(u)dV$$
where $dV$ is the volume element of $(M,g)$.

A map $u\in C^1(M,N)$ is called harmonic if it is a critical point of the energy $E$.

By Nash embedding theorem we know that $(N,h)$ can be isometrically into an Euclidean space $R^K$ with some positive integer $K$.
Then $(N,h)$ may be considered as a submanifold of $R^K$ with the metric induced from the Euclidean metric. Thus a map
$u\in C^1(M,N)$ can be considered as a map of $C^1(M,R^K)$ whose image lies on $N$. In this
sense we can get the following Euler-Lagrange equation
$$\triangle u=A(u)(du,du).$$

The tension field $\tau(u)$ is defined by
$$\tau(u)=\triangle_M u-A(u)(du,du)$$
where $A(u)(du,du)$ is the second fundamental form of $N$ in $R^K$.
So $u$ is harmonic means that $\tau(u)=0$.

The harmonic mappings are of special interest when $M$ is a Riemann surface.
Consider a sequence of mappings $u_n$ from Riemann surface $M$ to $N$ with bounded energies. It is clear that
$u_n$ converges weakly to $u$ in $W^{1,2}(M,N)$ for some $u\in W^{1,2}(M,N)$. But in general, it mayn't converge strongly
in $W^{1,2}(M,N)$. When $\tau(u_n)=0$, i.e. $u_n$ are all harmonic, Parker in \cite{parker} proved that the lost energy
is exactly the sum of some harmonic spheres which is defined as a harmonic mapping from $S^2$ to $N$. This result is called energy
identity. Also he proved that the images of these harmonic spheres and $u(M)$ are connected, i.e. there is no neck during
blowing up.

When $\tau(u_n)$ is bounded in $L^2$, the energy identity is
proved in \cite{Qing} for the sphere, in \cite{DT} and \cite{Wa}
for the general target manifold. In \cite{QT} they proved there is
no neck during blowing up. For the heat flow of harmonic mappings,
the results can also be found in \cite{To1,To2}. When the target
manifold is a sphere, in \cite{LZ} we proved the energy identity
for a sequence of mappings with tension fields bounded in $L\ln^+
L$ where they used good observations in \cite{LWa}. On the other
hang, in \cite{LZ} we constructed a sequence of mappings with
tension fields bounded in $L\ln^+ L$ such that there is positive
neck during blowing up. Furthermore, in \cite{zhu} the second
author proved the neckless during blowing up for a sequence of
maps $u_n$ with
$$\lim_{\delta\to 0}\sup_n \sup_{B(x,\delta)\subset D_1}\|\tau(u_n)\|_{L\ln^+ L(B(x,\delta))}=0.$$

In this paper we proved the energy identity and neckless during
blowing up of a sequence of maps $u_n$ with $\tau(u_n)$ bounded in
$L^p$ for some $p\geq \frac 65$, for the general target manifold.

When $\tau(u_n)$ is bounded in $L^p$ for some $p>1$, the small energy regularity proved in \cite{DT} implies that $u_n$
converges strongly in $W^{1,2}(M,N)$ outside a finite set of points. For simplicity in exposition, it is no matter to assume that
$M$ is the unit disk $D_1=D(0,1)$ and there is only one singular point at 0.

In this paper we proved the following theorem.
\begin{theorem}\label{bubble1}
Let $\{u_n\}$ be a sequence of mappings from $D_1$ to $N$ in $W^{1,2}(D_1,N)$ with tension field $\tau(u_n)$. If\\
(a) $\|u_n\|_{W^{1,2}(D_1)}+\|\tau(u_n)\|_{L^p(D_1)}\leq \Lambda$ for some $p\geq\frac 65$;\\
(b) $u_n\rightarrow u$ strongly in $W^{1,2}(D_1\setminus \{0\},R^K)$ as $n\rightarrow \infty$.

Then there exist a subsequence of $\{u_n\}$ (we still denote it by $\{u_n\}$) and some
nonnegative integer $k$. For any $i=1,...,k$, there exist points $x_n^i$, positive numbers $r_n^i$ and a nonconstant harmonic
sphere $w^i$ (which we view as a map from $R^2\cup\{\infty\}\rightarrow N$) such that\\
(1) $x_n^i\rightarrow 0,r_n^i\rightarrow 0$ as $n\rightarrow \infty$;\\
(2) $\lim_{n\rightarrow\infty}(\frac{r_n^i}{r_n^j}+\frac{r_n^j}{r_n^i}+\frac{|x_n^i-x_n^j|}{r_n^i+r_n^j})=\infty$ for any $i\neq j$;\\
(3) $w^i$ is the weak limit or strong limit of $u_n(x_n^i+r_n^i x)$ in $W^{1,2}_{Loc}(R^2,N)$;\\
(4) \textbf{Energy identity:}
\begin{equation}\label{energy identity}
\lim_{n\rightarrow \infty} E(u_n,D_1)=E(u,D_1)+\sum^k_{i=1}E(w^i);
\end{equation}
(5) \textbf{Neckless: }The image $u(D_1)\cup\bigcup^k_{i=1}w^i(R^2)$ is a connected set.
\end{theorem}

This paper is organized as follows. In section 2 we state some basic lemmas and some standard arguments in the blow-up analysis.

In section 3 and section 4 we prove theorem \ref{bubble1}. In the proof, we used delicate analysis on the difference between
normal energy and tangential energy. Energy identity is proved in section 3 and neckless is proved in section 4.

Throughout this paper, without illustration the letter $C$ denotes
a positive constant which depends only on $p,\Lambda$ and the
target manifold $N$ and may vary in different cases. Furthermore
we always don't distinguish the sequence and its subsequence.

\section{Some basic lemmas and standard arguments}

We recall the regular theory for the mapping with small energy on
the unit disk and the tension field in $L^p\hspace{2mm}(p>1)$.
\begin{lemma}\label{small}
Let $\bar{u}$ be the mean value of $u$ on the disk $D_{\frac 12}$. There exists a positive constant $\epsilon_N$ that depends only
on the target manifold such that if $E(u,D_1)\leq\epsilon^2_N$ then
\begin{equation}
\|u-\bar{u}\|_{W^{2,p}(D_{\frac12})}\leq C(\|\nabla u\|_{L^2(D_1)}+\|\tau(u)\|_p)\label{2p}
\end{equation}
where $p>1$.

As a direct consequence of (\ref{2p}) and the Sobolev embedding $W^{2,p}(R^2)\subset C^0(R^2)$, we have
\begin{equation}\label{OSC}
\|u\|_{Osc(D_{\frac12})}=\sup_{x,y\in D_{\frac12}}|u(x)-u(y)|\leq C(\|\nabla u\|_{L^2(D_1)}+\|\tau(u)\|_p).
\end{equation}
\end{lemma}

The lemma has been proved in \cite{DT}.

\textbf{Remark 1.} In \cite{DT} they proved this lemma for the mean value of $u$ on the unit disk. Note that
$$|\frac{\int_{D_1}u(x)dx}{|D_1|}-\frac{\int_{D_{\frac 12}}u(x)dx}{|D_{\frac 12}|}|\leq C\|\nabla u\|_{L^2(D_1)}.$$
So we can use the mean value of $u$ on $D_{\frac 12}$ in this lemma.

\textbf{Remark 2.} Suppose we have a sequence of mappings $u_n$ from the unit disk $D_1$ to $N$ with
$\|u_n\|_{W^{1,2}(D_1)}+\|\tau(u_n)\|_{L^p(D_1)}\leq \Lambda$ for some $p>1$.

A point $x\in D_1$ is called an energy concentration point (blow-up point) if for any $r,D(x,r)\subset D_1$,
$$\sup_n E(u_n,D(x,r))>\epsilon^2_N$$
where $\epsilon_N$ is given in this lemma.

If $x\in D_1$ isn't an energy concentration point, then we can find a positive number $\delta$ such that
$$E(u_n,D(x,\delta))\leq \epsilon^2_N, \forall n.$$
Then it follows from Lemma \ref{small} that we have a uniformly $W^{2,p}(D(x,\frac{\delta}{2}))$-bound for $u_n$.
Because $W^{2,p}$ is compactly embedded into $W^{1,2}$, there is a subsequence of $u_n$ (denoted
by $u_n$) and $u\in W^{2,p}(D(x,\frac{\delta}{2}))$ such that
$$\lim_{n\rightarrow\infty}u_n=u\textrm{ in }W^{1,2}(D(x,\frac{\delta}{2})).$$
So $u_n$ converges to $u$ strongly in $W^{1,2}(D_1)$ outside a finite set of points.\\

Under the assumptions in the theorems, by the standard blow-up argument, i.e. rescalling $u_n$ suitable and repeated, we can
obtain some nonnegative integer $k$. For any $i=1,...,k$, there exist a point $x_n^i$, a positive number $r_n^i$ and a
nonconstant harmonic sphere $w^i$ satisfying (1), (2) and (3) of the theorem 1. By the standard induction argument in
\cite{DT} we only need to prove the theorems in the case that there is only one bubble.

In this case we may assume that $w$ is the strong limit of the sequence
$u_n(x_n+r_nx)$ in $W^{1,2}_{Loc}(R^2)$. It does nothing to assume that $x_n=0$. Set $w_n(x)=u_n(r_nx)$.

As
$$\lim_{\delta\rightarrow 0}\lim_{n\rightarrow \infty}E(u_n,D_1\setminus D_\delta)=E(u,D_1),$$
 the energy identity is equivalent to
\begin{equation}\label{claim1}
\lim_{\delta\rightarrow 0}\lim_{n\rightarrow \infty}\lim_{R\rightarrow \infty}E(u_n,D_\delta\setminus D_{r_nR})=0.
\end{equation}

To prove the set $u(D_1)$ and $w(R^2\cup\infty)$ is connected, it is enough to show that
\begin{equation}\label{no neck}
\lim_{\delta\rightarrow 0}\lim_{n\rightarrow \infty}\lim_{R\rightarrow \infty}\sup_{x,y\in D_\delta\setminus D_{r_nR}}
|u_n(x)-u_n(y)|=0.
\end{equation}

\section{Energy identity}
\hspace{6mm}In this section, we prove the energy identity for the general target manifold when $p\geq \frac65$.

Assume that there is only one bubble $w$ which is the strong limit of $u_n(r_n\cdot)$ in $W^{1,2}_{Loc}(R^2)$. Let $\epsilon_N$ be
the constant in Lemma \ref{small}. Furthermore, by the standard argument of blow-up analysis we can assume that for any $n$,
\begin{equation}\label{maximal}
E(u_n,D_{r_n})=\sup_{r\leq r_n,D(x,r)\subseteq D_1}E(u_n,D(x,r))=\frac{\epsilon^2_N}{4}.
\end{equation}

By the argument in \cite {DT}, we can show
\begin{lemma}\label{tangential}(\cite{DT})
If $\tau(u_n)$ is bounded in $L^p$ for some $p>1$, then the tangential energy on the neck domain equals to zero, i.e.
\begin{equation}\label{tan}
\lim_{\delta\rightarrow 0}\lim_{R\rightarrow\infty}\lim_{n\rightarrow\infty}
\int_{D_\delta\setminus D_{r_nR}}|x|^{-2}|\partial_\theta u|^2dx=0.
\end{equation}
\end{lemma}
\textbf{Proof:} The proof is the same as that in \cite{DT}, we sketch it.

For any $\epsilon>0$, take $\delta,R$ such that for any $n$,
$$E(u,D_{4\delta})+E(w,R^2\setminus D_R)+\delta^{\frac{4(p-1)}p}<\epsilon^2.$$
It is no matter to suppose that $r_nR=2^{-j_n},\delta=2^{-j_0}$. When $n$ is big enough, for any $j_0\leq j\leq j_n$,
there holds (see \cite{DT})
$$E(u_n,D_{2^{1-j}}\setminus D_{2^{-j}})<\epsilon^2.$$

For any $j$, set $h_n(2^{-j})=\frac{1}{2\pi}\int_{S^1}u_n(2^{-j},\theta)d\theta$ and
$$h_n(t)=h_n(2^{-j})+(h_n(2^{1-j})-h_n(2^{-j}))\frac{\ln(2^jt)}{\ln 2},t\in [2^{-j},2^{1-j}].$$
It is easy to check that
$$\frac{d^2h_n(t)}{dt^2}+\frac 1t\frac{dh_n(t)}{dt}=0,t\in [2^{-j},2^{1-j}].$$

Consider $h_n(x)=h_n(|x|)$ as a map from $R^2$ to $R^K$, then $\triangle h_n=0$ in $R^2$. Set $P_j=D_{2^{1-j}}\setminus D_{2^{-j}}$ we have
\begin{equation}\label{ta1}
\triangle(u_n-h_n)=\triangle u_n-\triangle h_n=\triangle u_n=A(u_n)+\tau(u_n),x\in P_j.
\end{equation}

Taking the inner product of this equation with $u_n-h_n$ and integrating over $P_j$, we get that
$$\int_{P_j}|\nabla (u_n-h_n)|^2dx=-\int_{P_j}(u_n-h_n)(A(u_n)+\tau(u_n))dx+\int_{\partial P_j}(u_n-h_n)(u_n-h_n)_rds.$$
Note that by the definition, $h_n(2^{-j})$ is the mean value of $\{2^{-j}\}\times S^1$ and $(h_n)_r$ is independent of $\theta$.
So the integral of $(u_n-h_n)(h_n)_r$ on $\partial P_j$ vanishes.

When $j_0<j<j_n$, by Lemma \ref{small} we have
\begin{eqnarray*}
\|u_n-h_n\|_{C^0(P_j)}&\leq &\|u_n-h_n(2^{-j})\|_{C^0(P_j)}+\|u_n-h_n(2^{1-j})\|_{C^0(P_j)}\\
&\leq &2\|u_n\|_{Osc(P_j)}\\
&\leq &C(\|\nabla u_n\|_{L^2(P_{j-1}\cup P_j\cup P_{j+1})}+2^{\frac{2(1-p)j}{p}}\|\tau(u_n)\|_p)\\
&\leq &C(\epsilon+2^{-\frac{2(p-1)j}p})\\
&\leq & C(\epsilon+\delta^{\frac{2(p-1)}p})\leq C\epsilon.
\end{eqnarray*}

Summing $j$ for $j_0<j<j_n$, we have
\begin{eqnarray}
\lefteqn{\int_{D_{\delta}\setminus D_{2r_nR}}|\nabla (u_n-h_n)|^2dx= \sum_{j_0<j<j_n}\int_{P_j}|\nabla (u_n-h_n)|^2dx}\nonumber\\
&\leq& \sum_{j_0<j<j_n}\int_{P_j}|u_n-h_n|(|A(u_n)|+|\tau(u_n)|)dx\nonumber\\
&&+\sum_{j_0<j<j_n}\int_{\partial P_j}(u_n-h_n)(u_n-h_n)_rds\nonumber\\
&\leq&C\epsilon(\int_{D_{2\delta}\setminus D_{2r_nR}}(|\nabla u_n|^2+|\tau(u_n)|)dx
+\int_{\partial D_{2\delta}\cup \partial D_{2r_nR}}|\nabla u_n|ds)\nonumber\\
&\leq&C\epsilon(\int_{D_{2\delta}\setminus D_{2r_nR}}|\nabla u_n|^2dx+\delta^{\frac{2(p-1)}{p}}+\epsilon)\nonumber\\
&\leq &C\epsilon.\label{ta2}
\end{eqnarray}
Here we use the inequality $\int_{\partial D_{2\delta}\cup \partial D_{2r_nR}}|\nabla u_n|ds\leq C\epsilon$, which can be derived from
the Sobolev trace embedding theorem.

As $h_n(x)$ is independent of $\theta$, it can be shown that
$$\int_{D_{2\delta}\setminus D_{2r_nR}}|x|^{-2}|\partial_\theta  u_n|^2dx\leq
\int_{D_{2\delta}\setminus D_{2r_nR}}|\nabla (u_n-h_n)|^2dx\leq C\epsilon.$$
So this lemma is proved.

It is left to show that the normal energy on the neck domain also equals to zero. We need the following Pohozaev equality
which was first proved by Lin-Wang \cite{LW1}.
\begin{lemma}\label{pohozaev lemma}
(\textbf{Pohozaev equality}, \cite{LW1}, lemma 2.4, P374)\\
Let $u$ be a solution to
$$\triangle u+A(u)(du,du)=\tau(u),$$
then there holds
\begin{equation}\label{POHO}
\int_{\partial D_t}(|\partial_r u|^2-r^{-2}|\partial_\theta u|^2)ds=\frac 2t\int_{D_t}\tau\cdot(x\nabla u)dx.
\end{equation}
As a direct corollary, integrating it over $[0,\delta]$ we have
\begin{equation}\label{po1}
\int_{D_\delta}(|\partial_r u|^2-r^{-2}|\partial_\theta u|^2)dx=\int^\delta_0\frac 2t\int_{D_t}\tau\cdot(x\nabla u)dxdt.
\end{equation}
\end{lemma}
\textbf{Proof:} Multiplying both side of the equation by $x\nabla u$ and integrating it over $D_t$, we get
$$\int_{D_t}|\nabla u|^2dx-t\int_{\partial D_t}|\partial_r u|^2ds+\frac 12\int_{D_t}x\nabla|\nabla u|^2dx
=-\int_{D_t}\tau\cdot(x\nabla u)dx.$$
Note that
$$\frac 12\int_{D_t}x\nabla|\nabla u|^2dx=-\int_{D_t}|\nabla u|^2dx+\frac t2\int_{\partial D_t}|\nabla u|^2ds.$$
Hence,
$$\int_{\partial D_t}(|\partial_r u|^2-\frac 12|\nabla u|^2)ds=\frac 1t\int_{D_t}\tau\cdot(x\nabla u)dx.$$
As $|\nabla u|^2=|\partial_r u|^2+r^{-2}|\partial_\theta u|^2$, we proved this lemma.

Now we use this equality to estimate the normal energy on the neck domain. We prove the following lemma.
\begin{lemma}\label{normal}
If $\tau(u_n)$ is bounded in $L^p$ for some $p\geq \frac 65$, then for $\delta$ small enough, there holds
$$|\int_{D_\delta}(|\partial_r u_n|^2-|x|^{-2}|\partial_\theta u|^2)dx|\leq C\delta^{\frac{p-1}p}$$
where $C$ depends on $p$, $\Lambda$, the target manifold $N$ and the bubble $w$.
\end{lemma}
\textbf{Proof:} Take $\psi\in C_0^\infty (D_2)$ satisfying that $\psi=1$ in $D_1$, then
$$\triangle(\psi u_n)=\psi A(u_n)(du_n,du_n)+\psi\tau_n+2\nabla\psi\nabla u_n+u_n\triangle\psi.$$
Set $g_n=\psi A(u_n)(du_n,du_n)+\psi\tau_n+2\nabla\psi\nabla u_n+u_n\triangle\psi$. When $|x|<1$,
$$\partial_iu_n(x)=R_i\ast g_n(x)=\int \frac{x_i-y_i}{|x-y|^2}g_n(y)dy.$$

Let $\Phi_n$ be the Newtonian potential of $\psi\tau_n$, then $\triangle \Phi_n=\psi\tau_n$. The corresponding Pohozaev equality is
\begin{equation}\label{po2}
\int_{D_\delta}(|\partial_r \Phi_n|^2-r^{-2}|\partial_\theta \Phi_n|^2)dx=\int^\delta_0\frac 2t\int_{D_t}\psi\tau_n\cdot(x\nabla \Phi_n)dxdt.
\end{equation}
Here $\partial_i\Phi_n(x)=R_i\ast (\psi\tau_n)(x)=\int \frac{x_i-y_i}{|x-y|^2}(\psi\tau_n)(y)dy$.

As $\tau_n$ is bounded in $L^p\hspace{2mm}(p>1)$, there holds
$$\int_{D_\delta}|\nabla \Phi_n|^2dx\leq C\delta^{\frac{4(p-1)}p}\|\nabla \Phi_n\|^2_{\frac{2p}{2-p}}
\leq C\delta^{\frac{4(p-1)}p}\|\tau_n\|^2_p\leq C\delta^{\frac{4(p-1)}p}.$$

By (\ref{po2}), it can be shown that for any $\delta>0$,
\begin{equation}\label{po3}
|\int^\delta_0\frac 1t\int_{D_t}\psi\tau_n\cdot(x\nabla \Phi_n)dxdt|
\leq\int_{D_\delta}|\nabla \Phi_n|^2dx\leq C\delta^{\frac{4(p-1)}p}.
\end{equation}
For $\delta$ small enough, we have
\begin{eqnarray}
&&|\int_{D_\delta}(|\partial_r u_n|^2-r^{-2}|\partial_\theta u_n|^2)dx|\nonumber\\
&=&|\int^\delta_0\frac 2t\int_{D_t}\tau_n\cdot(x\nabla u_n)dxdt|\nonumber\\
&\leq&2|\int^\delta_0\frac 1t\int_{D_t}\tau_n\cdot(x\nabla \Phi_n)dxdt|
+2\int^\delta_0\frac 1t\int_{D_t}|x\tau_n||\nabla (u_n-\Phi_n)(x)|dxdt\nonumber\\
&\leq&C\delta^{\frac{4(p-1)}p}+2\int_{D_\delta}|x\tau_n||\nabla (u_n-\Phi_n)(x)|(\int^\delta_{|x|}\frac1tdt)dx\nonumber\\
&\leq&C\delta^{\frac{4(p-1)}p}+2\int_{D_\delta}|\tau_n||\nabla (u_n-\Phi_n)(x)||x|\ln \frac{1}{|x|}dx.\label{po4}
\end{eqnarray}

For any $j>0$, set $\varphi_j(x)=\psi(\frac{x}{2^{2-j}\delta})-\psi(\frac{x}{2^{-2-j}\delta})$. When $2^{-j}\delta\leq |x|<2^{1-j}\delta$,
we have
\begin{eqnarray}
|\partial_i(u_n-\Phi_n)(x)|&=&|\int \frac{x_i-y_i}{|x-y|^2}(g_n(y)-\psi\tau_n(y))dy|\nonumber\\
&\leq &\int \frac{|\psi A(u_n)(du_n,du_n)+2\nabla\psi\nabla u_n+u_n\triangle\psi|(y)}{|x-y|}\nonumber\\
&\leq &\int \frac{|\psi A(u_n)(y)|}{|x-y|}dy+C\int_{1<|y|<2}(|\nabla u_n|+|u_n|)(y)dy\nonumber\\
&\leq &\int \frac{|\varphi_j A(u_n)(y)|}{|x-y|}dy+\int \frac{|(\psi-\varphi_j)A(u_n)(y)|}{|x-y|}dy+C\nonumber\\
&\leq &\int \frac{|\varphi_j A(u_n)(y)|}{|x-y|}dy+\frac{\int |A(u_n)(y)|dy}{|x|}+C\nonumber\\
&\leq &\int \frac{|\varphi_j A(u_n)(y)|}{|x-y|}dy+\frac{C}{|x|}.\label{po5}
\end{eqnarray}

When $\delta>0$ is small enough and $n$ is big enough, for any $j>0$ we claim that
\begin{equation}\label{po6}
\|\varphi_j A(u_n)\|_{\frac{p}{2-p}}\leq C(2^{-j}\delta)^{-\frac{4(p-1)}{p}}
\end{equation}
where the constant $C$ depends only on $p$, $\Lambda$, the bubble $w$ and the target manifold $N$.

Take $\delta>0$ and $R(w)$ which depends on $w$ such that
$$E(u,D_{8\delta})\leq\frac{\epsilon^2_N}{8};\hspace{4mm}E(w,R^2\setminus D_{R(w)})\leq\frac{\epsilon^2_N}{8}.$$
The standard blow-up analysis (see \cite{DT}) show that for any $j$ with $8r_nR(w)\leq 2^{-j}\delta$ and $n$ big enough, there holds
$$E(u_n,D_{2^{4-j}\delta}\setminus D_{2^{-3-j}\delta})\leq\frac{\epsilon^2_N}{3}.$$
By (\ref{maximal}), when $2^{-j}\delta<\frac{r_n}{16}$, there holds
$$E(u_n,D_{2^{4-j}\delta}\setminus D_{2^{-3-j}\delta})\leq \frac{\epsilon^2_N}{4}.$$
So when $2^{-j}\delta<\frac{r_n}{16}$ or $2^{-j}\delta\geq 8r_nR(w)$, by Lemma \ref{small}, we have
\begin{eqnarray*}
\|\varphi_j A(u_n)\|_{\frac{p}{2-p}}&\leq &C\|\nabla u_n\|^2_{L^{\frac{2p}{2-p}}(D_{2^{3-j}\delta}\setminus D_{2^{-2-j}\delta})}\\
&\leq &C\|u_n-\overline{u_{n,j}}\|^2_{W^{2,p}(D_{2^{3-j}\delta}\setminus D_{2^{-2-j}\delta})}\\
&\leq &C[(2^{-j}\delta)^{-\frac{4(p-1)}{p}}\|\nabla u_n\|^2_{L^2(D_{2^{4-j}\delta}\setminus D_{2^{-4-j}\delta})}+\|\tau(u_n)\|^2_p]\\
&\leq &C(2^{-j}\delta)^{-\frac{4(p-1)}{p}}
\end{eqnarray*}
where $\overline{u_{n,j}}$ is the mean of $u_n$ on $D_{2^{3-j}\delta}\setminus D_{2^{-2-j}\delta}$.

On the other hand, when $\frac{r_n}{16}\leq 2^{-j}\delta\leq 8r_nR(w)$, we can find no more than $CR(w)^2$ balls
with radius $\frac{r_n}{2}$ to cover $D_{2^{3-j}\delta}\setminus D_{2^{-2-j}\delta}$, i.e.
$$D_{2^{3-j}\delta}\setminus D_{2^{-2-j}\delta}\subset \bigcup^m_{i=1}D(y_i,\frac{r_n}{2}).$$
Denote $B_i=D(y_i,\frac{r_n}{2})$ and $2B_i=D(y_i,r_n)$. By (\ref{maximal}), for any $i$ with $i\leq m$ there holds
$$E(u_n,2B_i)\leq \frac{\epsilon^2_N}{4}.$$
Using Lemma \ref{small} we have
\begin{eqnarray*}
\|\varphi_j A(u_n)\|_{\frac{p}{2-p}}&\leq &C\|\nabla u_n\|^2_{L^{\frac{2p}{2-p}}(D_{2^{3-j}\delta}\setminus D_{2^{-2-j}\delta})}\\
&\leq &C(\sum_{i=1}^m\|\nabla u_n\|^{\frac{2p}{2-p}}_{L^{\frac{2p}{2-p}}(B_i)})^{\frac{2-p}{p}}\\
&\leq &C\sum_{i=1}^m\|\nabla u_n\|^2_{L^{\frac{2p}{2-p}}(B_i)}\\
&\leq &C\sum_{i=1}^m\|u_n-\overline{u_{n,i}}\|^2_{W^{2,p}(B_i)}\\
&\leq &C\sum_{i=1}^m((r_n)^{-\frac{4(p-1)}{p}}\|\nabla u_n\|^2_{L^2(2B_i)}+\|\tau(u_n)\|^2_p)\\
&\leq &Cm((2^{-j}\delta)^{-\frac{4(p-1)}{p}}+1)\\
&\leq &C(2^{-j}\delta)^{-\frac{4(p-1)}{p}}
\end{eqnarray*}
where $\overline{u_{n,i}}$ is the mean of $u_n$ over $B_i$ and the constant $C$ depends only on $p$, $\Lambda$,
the bubble $w$ and the target manifold $N$. So we proved the claim (\ref{po6}).

By (\ref{po5}) and (\ref{po6}), when $p>1$ we get that
\begin{eqnarray}
&&\int_{D_\delta}|\tau_n||\nabla (u_n-\Phi_n)(x)||x|\ln \frac{1}{|x|}dx\nonumber\\
&\leq&\sum^{\infty}_{j=1}\int_{2^{-j}\delta<|x|<2^{1-j}\delta}|\tau_n||\nabla(u_n-\Phi_n)(x)||x|\ln \frac{1}{|x|}dx\nonumber\\
&\leq&C\sum^{\infty}_{j=1}\int_{2^{-j}\delta<|x|<2^{1-j}\delta}|\tau_n|(\frac{1}{|x|}+\int \frac{|\varphi_j A(u_n)(y)|}{|x-y|}dy)
|x|\ln \frac{1}{|x|}dx\nonumber\\
&\leq&C(\int_{D_\delta}|\tau_n|\ln \frac{1}{|x|}dx+
\sum^{\infty}_{j=1}\int_{2^{-j}\delta<|x|<2^{1-j}\delta}|\tau_n|(\int \frac{|\varphi_j A(u_n)(y)|}{|x-y|}dy)|x|\ln \frac{1}{|x|}dx)\nonumber\\
&\leq&C(\|\ln \frac{1}{|\cdot|}\|_{L^{\frac{p}{p-1}}(D_\delta)}
+\sum^{\infty}_{j=1}2^{-j}\delta\ln\frac{2^j}{\delta}\|\int \frac{|\varphi_j A(u_n)(y)|}{|\cdot-y|}dy\|_{\frac{p}{p-1}})\|\tau_n\|_p\nonumber\\
&\leq&C(\delta^2(\ln \frac 1\delta)^{\frac 1{p-1}}+
\sum^{\infty}_{j=1}2^{-j}\delta\ln\frac{2^j}{\delta}\|\varphi_j A(u_n)\|_{\frac{2p}{3p-2}}).\label{po7}
\end{eqnarray}
Here we use the fact that the fraction integral operator $I(f)=\frac{1}{|\cdot|}\ast f$ is bounded from $L^q(R^2)$ to
$L^{\frac{2q}{2-q}}(R^2)$ for $1<q<2$.

When $p\geq \frac 65$, i.e. $\frac{2p}{3p-2}\leq \frac{p}{2-p}$, by (\ref{po6}) there holds
\begin{equation}
\|\varphi_j A(u_n)\|_{\frac{2p}{3p-2}}\leq C(2^{-j}\delta)^{\frac{5p-6}{p}}\|\varphi_j A(u_n)\|_{\frac{p}{2-p}}
\leq C(2^{-j}\delta)^{\frac{5p-6}{p}-\frac{4(p-1)}{p}}\leq C(2^{-j}\delta)^{-\frac{2-p}{p}}.\label{po8}
\end{equation}
From (\ref{po7}) and (\ref{po8}) we have
\begin{eqnarray}
&&\int_{D_\delta}|\tau_n||\nabla (u_n-\Phi_n)(x)||x|\ln \frac{1}{|x|}dx\nonumber\\
&\leq&C(\delta^2(\ln \frac 1\delta)^{\frac 1{p-1}}
+\sum^{\infty}_{j=1}2^{-j}\delta\ln\frac{2^j}{\delta}\|\varphi_j A(u_n)\|_{\frac{2p}{3p-2}})\nonumber\\
&\leq&C(\delta+\sum^{\infty}_{j=1}2^{-j}\delta\ln\frac{2^j}{\delta}(2^{-j}\delta)^{-\frac{2-p}{p}})\nonumber\\
&\leq&C(\delta+\delta^{\frac{2(p-1)}{p}}\ln\frac 1\delta))\nonumber\\
&\leq&C\delta^{\frac{p-1}p}.\label{po9}
\end{eqnarray}
It is clear that (\ref{po4}) and (\ref{po9}) imply that
\begin{equation}\label{po10}
|\int_{D_\delta}(|\partial_r u_n|^2-r^{-2}|\partial_\theta u_n|^2)dx|\leq C\delta^{\frac{p-1}{p}}.
\end{equation}

Now we use these lemmas to prove the energy identity. Note that $w$ is harmonic, from Lemma \ref{pohozaev lemma} we see that
$\int_{D_R}(|\partial_r w|^2-r^{-2}|\partial_\theta w|^2)dx=0$ for any $R>0$. It is easy to see that
$$\lim_{R\rightarrow\infty}\lim_{n\rightarrow\infty}|\int_{D_{r_nR}}(|\partial_r u_n|^2-r^{-2}|\partial_\theta u_n|^2)dx|=
\lim_{R\rightarrow\infty}|\int_{D_R}(|\partial_r w|^2-r^{-2}|\partial_\theta w|^2)dx|=0.$$

Letting $\delta\rightarrow 0$ in (\ref{po10}), we obtain
\begin{eqnarray*}
&&\lim_{\delta\rightarrow 0}\lim_{R\rightarrow\infty}\lim_{n\rightarrow\infty}
|\int_{D_\delta\setminus D_{r_nR}}(|\partial_r u_n|^2-r^{-2}|\partial_\theta u_n|^2)dx|\\
&\leq&\lim_{\delta\rightarrow 0}\lim_{n\rightarrow\infty}
|\int_{D_\delta}(|\partial_r u_n|^2-r^{-2}|\partial_\theta u_n|^2)dx|+
\lim_{R\rightarrow\infty}\lim_{n\rightarrow\infty}
|\int_{D_{r_nR}}(|\partial_r u_n|^2-r^{-2}|\partial_\theta u_n|^2)dx|\\
&=&0.
\end{eqnarray*}
Using Lemma \ref{tangential} we obtain that the normal energy also vanishes on the neck domain. So the energy identity is proved.

\section{Neckless}
\hspace{6mm}In this section we use the method in \cite{QT} to prove the neckless during blowing up.

For any $\epsilon>0$, take $\delta,R$ such that
$$E(u,D_{4\delta})+E(w,R^2\setminus D_R)+\delta^{\frac{4(p-1)}p}<\epsilon^2.$$
Suppose $r_nR=2^{-j_n},\delta=2^{-j_0}$. When $n$ is big enough, the standard blow-up analysis show that for any $j_0\leq j\leq j_n$,
$$E(u_n,D_{2^{1-j}}\setminus D_{2^{-j}})<\epsilon^2.$$

For any $j_0<j<j_n$, set $L_j=\min\{j-j_0,j_n-j\}$. Now we estimate the norm $\|\nabla u_n\|_{L^2(P_j)}$. Denote
$P_{j,t}=D_{2^{t-j}}\setminus D_{2^{-t-j}}$ and take $h_{n,j,t}$ similar to $h_n$ in the last section but
$h_{n,j,t}(2^{\pm t-j})=\frac 1{2\pi}\int_{S^1}u_n(2^{\pm t-j},\theta)d\theta$. By an argument similar to the one used in deriving
(\ref{ta2}), we have, for any $0<t\leq L_j$,
\begin{equation}
\int_{P_{j,t}}r^{-2}|\partial_\theta  u_n|^2dx
\leq C\epsilon(\int_{P_{j,t}}|\nabla u_n|^2dx+(2^{t-j})^{\frac{2(p-1)}{p}})
+\int_{\partial P_{j,t}}|u_n-h_{n,j,t}||\nabla u_n|ds.\label{neckp1}
\end{equation}

Set $f_j(t)=\int_{P_{j,t}}|\nabla u_n|^2dx$, a simple computation shows that
$$f_j^\prime(t)=\ln 2(2^{t-j}\int_{\{2^{t-j}\}\times S^1}|\nabla u_n|^2ds+2^{-t-j}\int_{\{2^{-t-j}\}\times S^1}|\nabla u_n|^2ds).$$

As $h_{n,j,t}$ is independent of $\theta$ and $h_{n,j,t}$ is the mean value of $u_n$ at the two components of
$\partial P_{j,t}$, by Poincar\'{e} inequality we get
\begin{eqnarray*}
\int_{\partial P_{j,t}}|u_n-h_{n,j,t}||\nabla u_n|ds
&=&\int_{\{2^{t-j}\}\times S^1}|u_n-h_{n,j,t}||\nabla u_n|ds\\
&&+\int_{\{2^{-t-j}\}\times S^1}|u_n-h_{n,j,t}||\nabla u_n|ds\\
&\leq &(\int_{\{2^{t-j}\}\times S^1}|u_n-h_{n,j,t}|^2ds)^{\frac 12}(\int_{\{2^{t-j}\}\times S^1}|\nabla u_n|^2ds)^{\frac 12}\\
&&+(\int_{\{2^{-t-j}\}\times S^1}|u_n-h_{n,j,t}|^2ds)^{\frac 12}(\int_{\{2^{-t-j}\}\times S^1}|\nabla u_n|^2ds)^{\frac 12}\\
&\leq &C(2^{t-j}\int_{\{2^{t-j}\}\times S^1}|\nabla u_n|^2ds+2^{-t-j}\int_{\{2^{-t-j}\}\times S^1}|\nabla u_n|^2ds)\\
&\leq &Cf_j^\prime(t).
\end{eqnarray*}

On the other hand, by a similar argument as we did in obtaining (\ref{po10}), we have
\begin{equation}
|\int_{P_{j,t}}(|\partial_r u_n|^2-r^{-2}|\partial_\theta u_n|^2)dx|
\leq C((2^{t-j})^{\frac{p-1}{p}}+(2^{-t-j})^{\frac{p-1}{p}})\leq C(2^{t-j})^{\frac{p-1}{p}}.\label{neckp2}
\end{equation}

Since $|\nabla u|^2=|\partial_r u|^2+r^{-2}|\partial_\theta u|^2=2r^{-2}|\partial_\theta u|^2+(|\partial_r u|^2-r^{-2}|\partial_\theta u|^2)$,
by (\ref{neckp1}) and (\ref{neckp2}) we have
\begin{eqnarray*}
f_j(t)&\leq&2\int_{P_{j,t}}r^{-2}|\partial_\theta u_n|dx
+|\int_{P_{j,t}}(|\partial_r u_n|^2-r^{-2}|\partial_\theta u_n|^2)dx|\\
&\leq&C\epsilon(f_j(t)+(2^{t-j})^{\frac{2(p-1)}{p}})+Cf_j^\prime(t)+C(2^{t-j})^{\frac{p-1}{p}}\\
&\leq&C(\epsilon f_j(t)+2^{-\frac{(p-1)j}{p}}2^{\frac{(p-1)t}p}+f_j^\prime(t)).
\end{eqnarray*}
Take $\epsilon$ small enough and denote $\epsilon_p=\frac{p-1}{p}\ln2$, then for some positive constant $C$ big enough there holds
$$f_j^\prime(t)-\frac 1Cf_j(t)+Ce^{-\epsilon_p j}e^{\epsilon_p t}\geq 0.$$
It is no matter to assume that $\epsilon_p>\frac 1C$, then we have
$$(e^{-\frac tC}f_j(t))^\prime+Ce^{-\epsilon_p j}e^{(\epsilon_p-\frac 1C)t}\geq 0.$$
Integrating this inequality over $[2,L_j]$, we get
$$f_j(2)\leq C(e^{-\frac {L_j}C}f_j(L_j)+e^{-\epsilon_p j}\int^{L_j}_1e^{(\epsilon_p-\frac 1C)t}dt)
\leq C(e^{-\frac {L_j}C}f_j(L_j)+e^{-\epsilon_p j}e^{(\epsilon_p-\frac 1C)L_j}).$$
Note that $j\geq L_j$, there holds
$$f_j(2)\leq C(e^{-\frac {L_j}C}f_j(L_j)+e^{-\frac jC}).$$
Since the energy identity has been proved in the
last section, we can take $\delta$ small such that the energy on the neck domain is
less than $\epsilon^2$ which implies that $f_j(L_j)<\epsilon^2$. So we get
$$f_j(2)\leq C(e^{-\frac {L_j}C}\epsilon^2+e^{-\frac jC}).$$
Using Lemma \ref{small} on the domain $P_j=D_{2^{1-j}}\setminus D_{2^{-j}}$ when $j<j_n$, we obtain
$$\|u_n\|_{Osc(P_j)}\leq C(\|\nabla u_n\|_{L^2(P_{j-1}\cup P_j\cup P_{j+1})}+2^{\frac{2(1-p)j}p}\|\tau(u_n)\|_p)
\leq C(f_j(2)+e^{-2\epsilon_p j}).$$
Summing $j$ from $j_0$ to $j_n$ we obtain that
\begin{eqnarray*}
\|u_n\|_{Osc(D_\delta\setminus D_{2r_nR})}&\leq&\sum^{j_n}_{j=j_0}\|u_n\|_{Osc(P_j)}\\
&\leq& C\sum^{j_n}_{j=j_0}(f_j(2)+e^{-2\epsilon_p j})\\
&\leq& C\sum^{j_n}_{j=j_0}(e^{-\frac {L_j}C}\epsilon^2+e^{-\frac jC}+e^{-2\epsilon_p j})\\
&\leq& C(\sum^{\infty}_{i=0}e^{-\frac iC}\epsilon^2+\sum^{\infty}_{j=j_0}e^{-\frac jC})\\
&\leq& C(\epsilon^2+e^{-\frac {j_0}C})\\
&\leq& C(\epsilon^2+\delta^{\frac 1C}).
\end{eqnarray*}
Here we use the assumption that $\epsilon_p>\frac{1}{C}$. So we proved that there is no neck during the blowing up.

Acknowledgment The research was supported by NSFC 11071236,
11131007 and 11101372 and by PCSIRT.

\end{document}